# LA NOTION D'IRRATIONALITÉ SELON UN MATHÉMATICIEN DU $X^E$ SIÈCLE: ABŪ JAʿFAR AL-KHĀZIN

**Nicolas Farès**
*Équipe d'Étude et de Recherche sur la Tradition Scientifique Arabe* – CNRS - Liban et Université Libanaise
nfares55@hotmail.com

**Résumé**

Le *Traité* d'Abū Jaʿfar al-Khāzin ($X^e$ s.), intitulé "*Commentaire de l'introduction du dixième livre du traité d'Euclide*" ("*tafsīr ṣadr al-maqāla al-ʿāshira min kitāb Uqlīdis*") existe en huit manuscrits. Nous en présentons ici une étude basée sur une première édition que nous avons faite d'après des copies des manuscrits de Paris, de Leyde et de Tunis, Ahmadiya**.**

La lecture du texte d'al-Khāzin montre que ce mathématicien a effectué une étude profonde du livre X des *Éléments*, pour en rendre un aperçu global. Il en présente, en effet, un commentaire condensé dont l'originalité se fait sentir tant au niveau de la forme qu'à celui du fond. Nous allons essayer tant que possible de dégager ces originalités, en nous référant surtout au commentaire de Pappus du "livre X", aux travaux de R. Rashed sur l'histoire de l'algèbre, et à un travail récent de M. Ben Miled sur les commentaires arabes du livre X.

Il nous semble que le travail d'al- Khāzin concernant le livre X, se situe dans une tradition arabe identifiée, surtout grâce aux travaux de R. Rashed et caractérisée par:
- une lecture algébrique du travail géométrique d'Euclide,
- une interprétation numérique des notions euclidiennes, rendue possible grâce à l'introduction d'une "*unité*" pour chaque type de grandeurs.

**Mots clés:** un manuscrit mathématique du $X^e$ siècle; al-Khāzin; Lectures algébriques du livre X des *Éléments* d'Euclide; grandeurs et nombres irrationnels.

\*\*\*

**Summary**

The treatise of Abū Jaʿfar al-Khāzin ($X^{th}$ century), entitled "*Commentary on the introduction of the tenth book of the treatise of Euclid*" ("*tafsīr ṣadr al-maqāla al-ʿāshira min kitāb Uqlīdis*") exists in eight manuscripts. In this paper we present a study of this treatise, based on a first edition we have made, according to the copies of the manuscripts of Paris, Leiden and Tunis, Ahmadiya.

A reading of this treatise of al-Khāzin proves that this mathematician did a profound study of the $X^{th}$ book of Euclid *Elements*, in order to gain a global comprehension. In effect, he presented a condensed commentary, the originality of which can be felt both on the formal and the contents levels. The aim of our study is to try, as much as possible, to show this originality, referring specially to Pappus commentary of the "$X^{th}$ book", the works of R. Rashed on the history of algebra and a recent work of M. Ben Miled on the Arabic commentaries of the $X^{th}$ book.

It appears to us that the work of al-Khāzin concerning the $X^{th}$ book, is situated in the current of an Arabic tradition identified, mainly due to the works of R. Rashed, by:
- an algebraic representation of the geometrical work of Euclid,
- a numerical interpretation of the Euclidian notions, rendered possible by the introduction of a *unity* for every kind of magnitude.

**Keywords:** a mathematical manuscript of the $X^{th}$ century; al-Khāzin; algebraic lectures of the $X^{th}$ book of Euclid's *Elements*, irrational magnitudes and numbers.

\*\*\*

**1. Préface: sur le mathématicien al-Khāzin et le manuscrit de son Traité.**

L'étude la plus récente sur "le nom, la vie et les faits" d'al-Khāzin est donnée par R. Rashed dans son ouvrage *Les mathématiques infinitésimales du IX$^e$ au XI$^e$ siècle* (Rashed, 1996, vol. 1, pp. 737-739). Nous reprenons ici, brièvement quelques renseignements le concernant.

Abū Ja'far al-Khāzin est un éminent mathématicien et astronome qui a travaillé dans la première moitié du X$^e$ siècle et qui était toujours en vie vers les années soixante de ce siècle, 350 H./961. Auteur d'un travail important en théorie des nombres (Rashed, 1979, repris en 1984, pp. 195-255 et Anbouba, 1979) son nom a été évoqué ces dernières années à l'occasion de la résolution du *grand théorème de Fermat*[1]. En effet, il est bien connu que ce théorème, qui porte le nom du grand mathématicien français du XVII$^e$ siècle, a attendu jusqu'à 1994 pour être résolu; mais, peu de gens savent qu'il a été énoncé explicitement, dans les cas n = 3, et n = 4, par des mathématiciens du X$^e$ siècle; certains ont même essayé de le démontrer, comme al-Khāzin[2].

Le *Traité* de Abū Ja'far al-Khāzin, intitulé "*Commentaire de l'introduction du dixième livre du traité d'Euclide*" ("*tafsīr ṣadr al-maqāla al-'āshira min kitāb Uqlīdis*") existe en huit manuscrits: 1) Feyzullāh 1359/6 ; 2) Berlin 5924/3 ; 3) Paris 2467 (ff. 201 – 207) ; 4) Leyde, Or. 1024/6 (pp. 65 - 82) ; 5) Téhéran, Daniskada-i- Adab. ğ 284/3 ; 6) Tunis, Ahmadiya 5482/4 (67$^b$- 72$^a$) ; 7) Patna 2928/10 ; 8) Hayderabād, Asaf. Riyad. 331/5. (Sizgin, 1974, pp. 298-299).

Nous présentons ici une étude sommaire de ce *Traité*, basée sur une première édition du manuscrit que nous avons faite d'après les copies de Paris, de Leyde, et de Tunis.

**2. Introduction: importance du sujet du *Commentaire* d'al-Khāzin.**

L'importance du *Traité* d'al-Khāzin ne revient pas seulement à celle de son auteur mais, aussi, à celle du sujet qu'il traite. Les propos de R. Rashed à ce sujet sont significatifs: "*On ne comprend rien à l'histoire de l'algèbre si on ne soulignait les apports de deux courants de recherche qui se sont développés durant la période précédemment considérée* (autour du 10$^e$ s). *Le premier portait sur l'étude des quantités irrationnelles, soit à l'occasion d'une lecture du dixième livre des Éléments, soit, en quelques sorte, indépendamment*" (Rashed, 1997, vol. 2, p. 37).

---

[1] Les travaux d'al-Khāzin en théorie des nombres, ont surtout été évoqués à l'occasion de "la journée annuelle" de la Société Mathématique de France (juin, 1995), consacrée à la résolution du *grand théorème de Fermat* (Conférence de C. Houzel, "*Le Théorème de Fermat à travers l'histoire de l'analyse diophantienne*").

[2] R. Rashed consacre une bonne partie du 4$^e$ chapitre de son *Entre arithmétique et algèbre* à l'évocation des travaux de ce mathématicien en théorie des nombres. On y trouve notamment une édition de son traité sur les "triplets pythagoriciens" et une reproduction d'une démonstration (erronée), qui lui est attribuée, du "grand théorème de Fermat" dans le cas n=3. (Rashed. R, 1984, pp. 195-255). A. Anbouba souligne l'erreur commise par Wœpcke et admise sans discussion depuis, par les historiens des sciences, notamment par Sarton et Suter (voir aussi Youschkevitch, 1976 pp. 69 et 91) et, qui fait du mathématicien deux personnages distincts: 1) Abū Ja'far al-Khāzin ; 2) Abū Ja'far Muhammad Ibn al-Husayn; (Anbouba, 1978, note p. 99). Dans le même article Anbouba présente une biographie assez étendue du mathématicien. On y lit notamment que le mathématicien a son nom cité dans les œuvres biobibliographiques d'Ibn al-Nadīm et d'Ibn al-Qiftī, ainsi que par des mathématiciens éminents tels que Ibn 'Irāq (fin du X$^e$ siècle), al-Bīrunī (X$^e$ - XI$^e$ siècle), al-Khayyām (XI$^e$ - XII$^e$ siècle), Ibn Abī Shukr al-Maghribī (XII$^e$ siècle).

L'analyse du *Traité* d'al-Khāzin contribue à dessiner avec plus de précision le schéma du développement de la théorie des irrationnels, basée sur le livre X des *Éléments* d'Euclide, dans la tradition arabe. Elle contribue, notamment, à comprendre dans quelle mesure "*l'algèbre, inventée au neuvième siècle par les mathématiciens de Bagdad, a fait découvrir un nouveau type d'analogie entre les nombres (entiers) et les quantités continues*" (Houzel, 2002, p. 257).

Le livre X des *Éléments* étudie la commensurabilité et l'incommensurabilité des grandeurs (continues), la rationalité et l'irrationalité des segments et des aires. Il traite quelques classes d'irrationnels, *simples* ou *composés* par juxtaposition (adjonction) ou par retranchement, représentés par des segments de droite. Il considère l'aire du carré construit sur chacun de ces segments, ainsi que le côté du carré équivalent au rectangle dont les côtés sont pris dans l'ensemble de ces segments[3]. En termes modernes, les longueurs de ces segments et les aires de ces rectangles s'expriment, à partir de longueurs données, par l'addition, la soustraction, la multiplication et l'extraction de la racine carrée. Bien entendu, une telle interprétation est anachronique; car, "*depuis la découverte des grandeurs irrationnelles, la mathématique grecque antique séparait rigoureusement les nombres (entiers) des quantités continues que l'on rencontre en géométrie*" (Houzel, C. 2002. p. 257) et, le livre X, tel qu'il a été conçu par Euclide, est un chapitre de sa géométrie.

Le nombre élevé des commentaires du livre X à travers l'histoire semble être dû à la difficulté de ce livre[4]. Cette difficulté revient surtout à la lourdeur du langage mathématique utilisé[5].

"*Le commentaire du livre X d'Euclide*" est, à quelques modifications insignifiantes près, le titre commun à de nombreux traités, écrits par des mathématiciens arabes, entre le IX$^e$ et le XI$^e$ siècles[6]. Ces écrits traitent la notion de "grandeurs irrationnelles" dans un style souvent différent de celui d'Euclide.

En effet, le livre X traite les grandeurs irrationnelles en utilisant les techniques efficaces et fiables offertes par la géométrie. Or, les commentateurs arabes de ce livre disposaient d'une nouvelle discipline mathématique, l'algèbre, accueillie et utilisée sans réserve par les contemporains et les successeurs directs d'al-Khwārizmī[7]. Ceux là, armés des opérations (théorie des équations et calculs) effectuées sur la "chose" quelle que soit la nature mathématique (arithmétique ou géométrique) de cette chose, ont appliqué sans aucune gêne l'algèbre à la théorie des grandeurs irrationnelles, ainsi qu'à divers autres domaines.

---

[3] Cf. (Ben Miled, 1999) où l'auteur donne un résumé de ce livre, condensé et «neutre», évitant toute interprétation. Voir aussi (Ben Miled, 2005, pp. 13-25).

[4] J. Itard écrit : "*Sa lecture demande au mathématicien moderne une solide préparation et un courage assuré*" ; (Itard, J. 1984, p. 94). Simon Stevin (1548-1620) avait déjà exprimé la même caractéristique du livre X en ces termes: "La difficulté du dixiesme livre d'Euclide est à plusieurs devenue en horreur, voire jusque à l'appeler la croix des mathématiciens, matière trop dure à digérer, et en laquelle ils n'aperçoivent aucune utilité" (Heath, 1956, vol. 3, p. 9).

[5] Un exemple de cette lourdeur est offert par la proposition X, 17, dont l'énoncé est le suivant: "*Si deux droites sont inégales, et qu'un parallélogramme égal à la quatrième partie du carré sur la plus petite soit appliqué sur la plus grande par défaut d'une figure carrée et qu'il la divise en segments commensurables en longueur, la plus grande sera, en puissance, plus grande que la plus petite par un carré sur une droite commensurable en longueur avec elle-même. Et…* " (énoncé semblable, de plusieurs lignes, pour la réciproque). Cf. B. Vitrac, *Euclide d'Alexandrie* (Vitrac, 1998, p. 142). Voici une paraphrase de l'énoncé et de sa réciproque: Soient $a$ et $b$ deux segments tels que $a > b$. On prend sur $a$ un segment $x$ de telle sorte que $x.(a-x) = ¼ b^2$, alors $x$ est commensurable avec $a-x$, si et seulement si, $a$ est commensurable avec le côté du carré équivalent à $a^2 - b^2$.

[6] Une liste des commentateurs (connus de nos jours) a été donnée par Marouane Ben Miled dans son article: "Les commentaires d'Al-Māhānī " (Ben Miled, 1999). Voir aussi (Rashed, 1997, vol. 2, pp. 37-38) et (Ben Miled, 2005, pp. 33-35).

[7] Cf. (Rashed, 1997, vol. 2, p. 34).

Les recherches sur l'histoire de l'algèbre, effectuées dans les dernières décennies[8], ont permis de reconnaître le rôle joué, dans la tradition arabe, par l'application de l'algèbre à la théorie euclidienne des irrationnels:
1) dans la confirmation de l'efficacité des outils offerts par l'algèbre, donc dans l'enracinement et le développement de l'algèbre elle-même,
2) dans le développement de la théorie des irrationnels d'Euclide, dont l'étude n'appartient plus seulement au domaine de la géométrie; ses principaux éléments étant désormais des racines d'équations algébriques.

L'application de l'algèbre à la théorie des irrationnels du livre X a ainsi donné un nouveau souffle à la théorie des irrationnels (algébriques); elle a été favorisée par l'interprétation numérique des irrationnels, rendue possible grâce au choix d'une longueur "*unité*"; elle a, à son tour, favorisé une telle interprétation[9].

**3. Remarques sur le contenu du Traité**.

Le *Traité* n'apporte presque pas de renseignement historique. A part le nom d'Euclide, al-Khāzin ne cite que celui de Sulaymān Ibn ʻIṣma, (mathématicien du X$^e$ siècle) qui avait travaillé dans le même domaine. Ce nom, sans d'autre indication sur le travail mentionné, est cité juste à la fin du Traité.

La lecture du *Traité* d'al-Khāzin montre qu'il ne s'agit pas d'une étude systématique des propositions du livre X. Le but du *Traité* est, comme l'indique le titre, d'expliquer l'introduction du livre X. Il insiste ainsi sur les notions de base et les définitions, mais aussi sur quelques résultats jugés essentiels par l'auteur. On a donc affaire à un résumé du livre X, condensé et original. Citons quelques indices de cette originalité, au niveau de la forme, qui reflètent des différences au niveau du fond avec le livre X.
1) Le non respect par al-Khāzin, de l'ordre dans lequel Euclide présente les notions et les propositions, est immédiatement perceptible. Au début du *Traité*, en effet, al-Khāzin présente les notions de rationalité (resp. d'irrationalité) avant celles de la commensurabilité (resp. d'incommensurabilité), contrairement à l'ordre suivi par Euclide. Plus loin, il présente les segments binomiaux, avant les bimédiaux. La dernière proposition du Livre X (prop. 115) est l'une des premières propositions du *Traité*; al-Khāzin la présente, dans l'un des premiers paragraphes[10], semble-t-il, pour illustrer sa conception de ce qu'il appelle "*rangs*" des irrationnels simples.
2) L'utilisation d'une terminologie étendue, distincte de celle d'Euclide, est une des propriétés du *Traité*. Nous en citons, non dans leur ordre d'importance, les termes: grandeurs "*ʻiẓam*" (pl. "*aʻẓām*"), *genre* ("*ğins*"), *valeur* ("*miqdār*"), *quantité continue*, grandeurs *homogènes*, *mesure* ou *évaluation* ("*qadr*", "*ʻiyār*", "*taqdīr*"), *unité* de valeur, *attribution* ("*iḍāfa*"), *produit*, *somme*, *racine* (*carrée*), *valeurs rationnelles* (resp. *irrationnelles*), segments rationnels *en longueur* et *en puissance*, *rangs* (d'irrationnels)... Certains termes appartiennent à des domaines mathématiques dépassant le cadre du livre X. Les énoncés de plusieurs propositions du Livre X sont adaptés à cette terminologie. Ainsi, dans l'énoncé de la prop. X, 21, al-Khāzin rend l'expression d'Euclide "*la ligne qui peut produire*" une aire par "*la racine du carré*" qui est égal à cette aire[11]. Un autre exemple plus significatif est celui où al-Khāzin parle de la "*racine du premier binomial*"[12]; l'énoncé de la prop. X, 54: "*Si une aire est contenue par une droite exprimable et par une binomiale première, la droite pouvant produire cette aire est irrationnelle, celle appelée binomiale*"[13] se trouve transformé par al-Khāzin en l'énoncé suivant: "*la racine du premier binomial est un binomial*".

---

[8] Cf. (Ahmad et Rashed, 1972), (Rashed, 1984), (Rashed 1997, vol. 2, ch. 2), etc… L'étude la plus récente de l'influence de l'algèbre arabe sur le chapitre concernant les irrationnels est, à notre connaissance présentée par Marouane Ben Miled, à l'occasion de l'édition du commentaire d'al-Māhānī (IX$^e$ siècle) et d'un commentaire anonyme du livre X, (Ben Miled, 1999); voir aussi (Ben Miled, 2005, pp. 35-92). S'il se confirme que ce commentaire anonyme revient à al-Māhānī (ce qui est bien probable d'après plusieurs arguments présentés par B. Miled), cela prouvera que même dès les premiers débuts de l'algèbre (précisément dès la deuxième moitié du IX$^e$ siècle), cette science a constitué un outil très efficace au développement de la théorie des irrationnels.
[9] Cf. (Houzel, 2002, pp. 257-259).
[10] Cf. le texte du *Traité* (Paris, f. 202r; Leyde, p. 67).

D'autres termes nouveaux sont introduits pour indiquer des notions qui ne se trouvent pas explicitement chez Euclide; c'est le cas de la notion de *rangs* d'irrationnels "*simples*"[14], absente du livre X. Après avoir affirmé qu'un certain segment de droite "*est irrationnel en longueur et en puissance, ainsi que ce qui le suit jusqu'à l'infini*", al-Khāzin dit que la catégorie des rationnels en puissance est le "*premier rang*" d'irrationnels et que, pour Euclide, les rangs d'irrationnels sont tous de la troisième catégorie (c'est-à-dire, des irrationnels dans l'absolu: en longueur et en puissance). Al-Khāzin parle de "rangs" au pluriel; d'après le contexte, il sous-entend qu'il y en a une infinité. Pourtant, al-Khāzin ne va pas plus loin dans la classification des irrationnels; dans la suite du *Traité*, il n'évoquera plus les irrationnels de rangs supérieurs à deux. Mais, quand il rappelle que, pour Euclide "*les rangs des sourds sont tous de la troisième catégorie, appelée médiale*,...," il sous-entend, semble-t-il, qu'Euclide a, pour une raison ou une autre, négligé (ou peut-être ignoré) les autres rangs[15] d'irrationnels simples.

3) La généralisation de la notion de rationalité à toutes les grandeurs est un des soucis du *Traité*. L'expression "*quantité continue*" ne se trouve pas explicitement dans le livre X. On la retrouve dans le commentaire de Pappus où la notion de continuité se trouve profondément expliquée et commentée (Thomson, 1930, p. 193). Une "*grandeur*" ("*'iẓam*") est, pour al-Khāzin, une "*quantité continue*"; elle possède une certaine "*valeur*" ("*miqdār*"). Or, la valeur d'une grandeur est exprimée en fonction d'une grandeur donnée qui lui est "*homogène*"; d'où la nécessité d'une unité de valeur pour chaque type de grandeur. La notion d'unité de valeur est bien explicite dans le texte d'al-Khāzin. Elle est appelée "*mesure*" par al-Khāzin[16] qui, plus loin[17], l'appelle, carrément, *l'unité* (*wāḥid* ou *al-wāḥid*).

On peut, ainsi, percevoir la visée généralisatrice d'al-Khāzin qui, en expliquant la notion de rationalité, ne voulait pas se limiter aux deux grandeurs: "longueur" et "surface", du livre X. De là on peut comprendre sa définition de la rationalité qui est unique pour toutes les grandeurs, ce qui n'est pas le cas dans le livre X. En effet, Euclide ne définit pas les segments rationnels et les surfaces rationnelles de la même façon. En fait, pour Euclide, la définition des surfaces rationnelles est la même que celle utilisée de nos jours, tandis que la définition des segments rationnels est plus large que celle que nous utilisons de nos jours[18]. Quant à la rationalité au sens d'al-Khāzin, elle est la même, pour toutes les grandeurs, y compris les longueurs et les surfaces[19] et, c'est la nôtre[20].

---

[11] cf. le texte du *Traité* (Paris, f. 202v; Leyde p. 69).

[12] cf. le texte du *Traité* (Paris, f. 205r; Leyde p. 79).

[13] Cf. (Vitrac, 1998, vol. 3, p. 251).

[14] Ou monômes, c'est-à-dire, qui (dans un langage moderne) peuvent être représentés par *a.u*, $\sqrt{a}\,u, \sqrt{\sqrt{a}}\,u$, ... où *a* est un nombre rationnel positif ($\sqrt{a} \notin \mathbb{Q}$) et *u* est l'unité de longueur.

[15] En fait, la proposition 115 du livre X ne dit pas explicitement que *BA, AT, TM*, ..., sont de *rangs distincts* mais, qu'ils sont *distincts*, tout court. Le terme "rang" n'est pas explicite chez Pappus, bien que la notion existe; il n'est pas explicite non plus chez al-Māhānī ; celui-ci va pourtant bien plus loin en introduisant d'irrationnels absents chez Euclide et en utilisant (implicitement) leurs rangs; cf. (Ben Miled, 1999, p. 142 et 2005, p. 287).

[16] Les termes arabes correspondants utilisés par al-Khāzin ("*mi'yār*" ou "*mikyāl*"), que nous avons traduits tous deux par "*mesure*" peuvent être considérés comme désignant une unité de mesure, au sens général du terme (qu'il s'agisse d'un volume, d'une longueur, d'un poids, etc...). Cf. le texte du *Traité*, (Paris, f. 201v; Leyde, p. 66).

[17] Cf. le texte du *Traité* (Paris, f. 201v; Leyde, p. 66).

Cette différence avec Euclide, n'influe pas essentiellement sur la suite du travail d'al-Khāzin. Elle rend, pourtant, la terminologie du livre X plus cohérente et facilite sa lecture que le double usage du terme "rationnel" alourdit davantage. Ce n'est que plus loin qu'al-Khāzin introduit les deux types de rationalité pour les lignes (i.e. pour les segments de droite): *rationalité en longueur* et *rationalité en puissance*[21] et, cela conformément à la terminologie d'Euclide concernant la commensurabilité des segments[22].

4) Contrairement au style du livre X, les applications et exemples numériques sont présents partout dans le texte. L'interprétation numérique des notions euclidiennes, grâce à l'introduction de "*l'unité*", a permis à al-Khāzin de traduire en termes de l'algèbre plusieurs notions du livre X et de construire des irrationnels[23] ou de caractériser des binomiaux, en utilisant les nombres entiers et leurs fractions[24]; on est évidemment, bien loin de la notion de nombres rationnels mais, il s'agit d'un usage confondant intuitivement la *grandeur* mesurée avec sa valeur par rapport à l'unité de mesure.

Dans une des proposition, il présente un algorithme itératif pour calculer le terme d'une certaine suite en fonction du terme précédent[25].

---

[18] Si (dans un langage moderne) nous appelons $u$ l'unité de longueur et $u^2$ l'unité de surface, pour Euclide, une surface $S$ est rationnelle, si et seulement si, $S = a.u^2$ où $a \in \mathbb{Q}$ (plus précisément, $a \in \mathbb{Q}^+$); pourtant, une longueur $l$ est rationnelle, si et seulement si, $l^2 = a.u^2$ où $a \in \mathbb{Q}$ (par exemple, $\sqrt{2}$ (i.e. $\sqrt{2}u$) est rationnel dans la terminologie d'Euclide, cf. Vitrac, B. 1998, *Définitions*, p. 101, dé. 3).

[19] Cf. le texte du *Traité* (Paris, f. 201v, li: 4-10; Leyde p. 65, li: 7-18).

[20] Al-Khāzin est explicite sur le fait qu'un segment rationnel en puissance (tel que $\sqrt{2}u$) est *sourd*: "*il est le premier rang des sourds*" (Paris, f. 202v, li: 10; Leyde p. 69, li: 3).

[21] Cette terminologie existe déjà chez Pappus; cf. (Thomson, 1930, p. 195).

[22] Cf. le texte du *Traité* (Paris, f. 201v, 202r, 202v; Leyde, p. 66 et 68).

[23] Cf. le texte du *Traité* (Paris, f. 206v et suite; Leyde, p. 88 et 89).

[24] Cf. le texte du *Traité* (Paris, f. 206v; Leyde pp. 88-89).

[25] Cf. le texte du *Traité*, (Paris, ff. 202r-202v; Leyde, pp. 67-69; Tunis, ff. 66r-67r). Pour montrer que les *rangs d'irrationnels simples* sont en nombre infini, al-Khāzin considère la proposition X. 115, qu'il traite en faisant intervenir son style numérique: il prend un segment de droite $AB$ ($= u = u_0$) et considère un segment qui lui est commensurable $BC = b = 2u_0$; puis, en supposant que $AB$ est le segment unité, il construit une suite $(u_i)_{i \geq 1}$, infinie (dénombrable), de segments qui lui sont incommensurables (en longueur), dont chacun est incommensurable avec le précédent, et vérifiant les relations $u_{n+1}^2 = b\,u_n\,(= 2u_0 u_n)$, pour tout $n \in \mathbb{N}^*$:

$$(u_0 = u = 1;\ u_0 = u = 1\ ;\ u_1 = \sqrt{2}u_0\ ;\ u_2 = \sqrt{2\sqrt{2}}\,u_0\ ;\ u_3 = \sqrt{2\sqrt{2\sqrt{2}}}\,u_0\ ; \ldots).$$

Le commentaire d'al-Khāzin se présente ainsi sous un aspect numérique distinct de celui du livre X, qui est géométrique.

5) Le langage algébrique est utilisé tout au long du *Traité*. Il en est de même pour les opérations algébriques (addition, soustraction, multiplication, extraction de la racine carrée) appliquées aux nombres, ainsi qu'aux segments et aux aires. Après avoir multiplié des segments de droite, al-Khāzin utilise le "produit de deux rectangles". Il parle explicitement de *rectangle*s et de *carré*s qui sont des binomiaux, puis de "*la racine carrée d'un binomial*"[26]. Cela n'est pas admis dans le style du livre X qui exige le respect de l'homogénéité, et qui interdit la confusion entre des notions introduites pour des longueurs et d'autres introduites pour des surfaces.

### 4. Conclusion: situation mathématique du *Traité* d'al-Khāzin.

Le schéma que nous venons de présenter permet de situer la mathématique du *Traité* d'al-Khāzin avec plus de précision. Il ne s'agit pas d'un "commentaire" du livre X, comme l'indique son titre, mais d'une interprétation des notions les plus importantes de ce livre, dans un style influencé par une tradition en développement depuis la fondation de l'algèbre avec al-Khwārizmī. Ce style se caractérise par l'aspect numérique dans lequel sont présentés les notions et les exemples et, par l'application de l'algèbre à toutes les branches mathématiques; il se distingue nettement du style euclidien, qui est géométrique. Deux remarques peuvent, pourtant, être avancées et permettent de constater un certain recul de ce traité par rapport à la tradition algébrique susmentionnée:

1) Al- Khāzin n'a pratiquement pas évoqué d'autres irrationnels que ceux d'Euclide[27], se limitant vraisemblablement à la tâche exprimée par le titre du *Traité*, à savoir "expliquer l'essentiel" du livre X. Il est en effet peu probable qu'al-Khāzin ignorât l'existence d'irrationnels d'autres types, bien connus au X$^e$ siècle: il a lui-même parlé de rangs d'irrationnels; de plus, l'étude d'al- Māhānī (milieu du 9$^e$ s.), bien antérieure à celle d'al-Khāzin, prenait clairement en considération une multitude d'autres types d'irrationnels, simples ou composés[28] et, Pappus lui- même dans son commentaire du livre X qui est connu dans la tradition arabe depuis le premier quart du X$^e$ siècle, évoquait d'autres irrationnels, (trinomiaux ou polynomiaux) qu'Apollonius aurait étudiés[29]. Près d'un siècle après al-Khāzin, al-Karajī (début du 11$^e$ s.) est explicite sur l'insuffisance des irrationnels d'Euclide; il évoque les polynomiaux ainsi que d'autres types d'irrationnels (algébriques)[30]: "*je dis que les monômes sont infinis: le premier est rationnel absolument, comme cinq, le deuxième est le rationnel en puissance, comme racine de dix, le troisième est défini par rapport à son cube,*

---

Le calcul des $u_i$ fait intervenir une suite intermédiaire $(s_i)_{i \geq 1}$, d'aires rectangulaires: $s_i = bu_i \ (=2u_0 u_i) = (u_{i+1})^2$. Au cours de son calcul de $u_3$, al-Khāzin s'arrête pour déduire une formule de récurrence permettant de calculer $s_{n+1}$ en fonction de $s_n$; cette formule peut s'écrire de la suivante :

$$s_n = (u_{n+1})^2 = b \cdot \sqrt{s_{n-1}} = \sqrt{b^2 \cdot s_{n-1}} = \left( b^{2^n} \cdot s_{n-1}^{2^{n-1}} \right)^{\frac{1}{2^n}}.$$

[26] Cf. le texte du *Traité* (Paris, f. 205r; Leyde, p. 78).

[27] A part ceux qui s'expriment sous la forme $m^{\frac{1}{2^n}}$ où *m* est un rationnel et *n* est un entier.

[28] Cf. (Ben Miled, 1999, p. 142 et 2005, p. 287).

[29] Cf. (Thomson, 1930, pp. 192, 211, 219).

[30] Cf. (Ahmad et Rashed, 1972, p. 50).

*comme le côté de vingt, le quatrième est la médiale, définie par rapport à son carré carré, comme la racine de la racine de dix, le cinquième est le côté du quadrato-cube, ensuite le côté du cubo-cube, et [ces monômes] se divisent ainsi à l'infini"*[31].

2) Al-Khāzin utilise bien la terminologie algébrique traditionnelle de l'époque ainsi que les opérations de l'algèbre appliquées aux binomiaux. Il se soucie peu de l'homogénéité et de la nature des grandeurs irrationnelles considérées, multipliant deux segments mais aussi, deux rectangles, ce qui constitue un progrès sur le chemin de l'algébrisation de la théorie. Mais, l'algèbre, en tant que science qui a pour objet de poser et de résoudre les équations algébriques, est tout à fait absente du traité, même là où l'on croit qu'il serait naturel de la rencontrer. Par exemple, en traitant le problème de la division d'un segment "*AE* en deux parties dont le produit de l'une par l'autre est le quart du carré de *EC*"[32], il ne la traite pas dans le style des algébristes de l'époque[33]. Ce fait est à remarquer car, la traduction algébrique du problème (en l'équation: $x^2 + c = bx$), a été utilisée par son contemporain al-Ahwāzī[34], en même temps que la méthode géométrique d'Euclide. De plus, la traduction en équations du second degré des problèmes posés par l'extraction des racines carrées des apotomes (donc, probablement, des autres binomiaux) et leurs résolutions par les algorithmes d'al-Kwārizmī, avaient déjà été pratiquées d'une façon systématique, plus d'un demi-siècle plus tôt, par al-Māhānī[35]. En outre, al-Khāzin lui-même devait être familier de la théorie des équations algébriques car, en plus de ses travaux mathématiques en théorie des nombres et en géométrie, on lui doit une contribution en algèbre, remarquée par al-Khayyām, à savoir la résolution d'une équation cubique posée par al-Māhānī[36]. Al-Khayyām précise que cette résolution a été faite au moyen de sections coniques, ce qui la classe dans le courant de l'algèbre "géométrique".

Le *Traité* d'al-Khāzin appartient à la longue tradition des commentaires arabes du Livre X. Celle-ci est, en général, marquée par la traduction algébrique des notions de ce livre[37], associant d'une façon spontanée et intuitive la "grandeur" rationnelle *mu* et le "nombre" (entier) *m* (*u* étant l'unité de grandeur convenable), ainsi que la fraction rationnelle $\frac{mu}{n}$ et le rapport $\frac{m}{n}$ des deux nombres entiers *m* et *n*. C'est ainsi que, vers la fin du *Traité*, al-Khāzin construit, à l'aide de "nombres" (entiers) des irrationnels euclidiens[38]. Plus tard, al-

---

[31] Cf. *Al-Bāhir en Algèbre* (Ahmad et Rashed, 1972, p. 40), et *al-Bad'fī al-hisāb* (Al-Karajī, p. 29).

[32] Cf. le texte du *Traité* (Paris, f. 205v; Leyde p. 80).

[33] En posant une de ces parties, par exemple, TE comme "chose" et l'autre, AT, comme AE moins cette "chose" pour transformer le problème en une équation du second degré et le résoudre en appliquant l'algorithme convenable d'al-Kwārizmī. En fait, al-Khāzin utilise deux méthodes pour résoudre ce problème : La première, qu'il qualifie de géométrique, est, pratiquement, celle d'Euclide (proposition X, 54). Dans la deuxième, qu'il appelle "*méthode par le calcul*", il ne se réfère pas à l'algorithme d'al-Khwārizmî mais, à la proposition II, 5, d'Euclide.

[34] Cf. (Al-Ahwāzī, Paris, Or. 2467, f. 209r).

[35] Cf. (Ben Miled, 2005, p. 297-333).

[36] Il s'agit d'une équation du troisième degré à laquelle al-Māhānī rend un problème du *De la sphère et du cylindre* d'Archimède. Al-Khayyām dit qu'il a fallu attendre, jusqu'à ce que abū-Ja'far al-Khāzin a eu le génie de la résoudre par les sections coniques; voir (Rashed et Vahabzadeh, 1999, texte arabe, p. 117, trad. fr. p. 116).

[37] Une des exceptions à cette règle fût, d'après R. Rashed (Ahmad et Rashed, 1972, p. 40), le commentaire d'Ibn al-Haytham.

[38] "*Ainsi, il s'est montré, de ce que nous avons décrit que, pour tout nombre, si tu retranches de son carré le quart, alors le nombre avec la racine du reste est un premier binomial; et, pour tout nombre, si tu ajoutes à son carré le tiers, alors le nombre, avec la racine de la somme, est un second binomial; et, de*

Karajī, essaye en quelque sorte, de légaliser cette pratique[39] qui date depuis les premiers commentaires dont celui d'al-Māhānī. Près d'un siècle après al-Karajī, on remarque chez al-Khayyām une nette approche entre la notion de rapport et celle de nombre[40], qui pourrait marquer le début de l'idée de "nombre rationnel". Tout ce mouvement causé par l'introduction de l'algèbre au 9$^e$ siècle et de son application à tous les domaines mathématiques, constitue un progrès remarquable vers la reconnaissance et l'introduction des nombres rationnels et des nombres irrationnels (algébriques).

---

*tout nombre, si tu retranches le quart, la racine du nombre avec la racine du reste est un troisième binomial; et, pour tout nombre, si tu retranches de son carré la moitié, alors le nombre et la racine du reste est un quatrième binomial; et, pour tout nombre, si tu doubles le carré, alors le nombre avec la racine du total est un cinquième binomial; et si, de tout nombre tu en retranches la moitié, la racine du nombre et la racine du reste est un sixième binomial*" (cf. le texte du *Traité* : Paris, f. 206v, ou Leyde, pp. 88-89):

[39] Al-Karajī exprime explicitement son projet en ces termes: "*et moi, je vous montre comment transposer ces dénominations* (celles des irrationnels d'Euclide) *aux nombres et j'en ajoute parce qu'ils ne sont pas suffisants pour le calcul*": "wa anā urīka naqla hādhihi al-alqāb ila al ʿadad wa azīdu fīhā li'annahu lā yuktafa bihā fī al-hisāb".

[40] Cf. (Rashed et Vahabzadeh, 1999, p. 276 et p. 341, lignes 10-20).